\documentclass[11pt]{amsart}
\usepackage{amsfonts,amssymb,amscd,amsmath,enumerate,verbatim,calc,latexsym,pstcol,pst-plot,pst-3d,}

\def\NZQ{\mathbb}               

\def\QQ{{\NZQ Q}}
\def\ZZ{{\NZQ Z}}
\def\RR{{\NZQ R}}

\def\frk{\mathfrak}               

\def\Phi{{\frk N}}
\def\a{{\alpha}}
\def\b{{\beta}}

\def\d{{\delta}}
\def\opn#1#2{\def#1{\operatorname{#2}}} 
\opn\chara{char} \opn\length{\ell} \opn\pd{pd} \opn\rk{rk}
\opn\projdim{proj\,dim} \opn\injdim{inj\,dim} \opn\rank{rank}
\opn\depth{depth} \opn\grade{grade} \opn\height{height}
\opn\embdim{emb\,dim} \opn\codim{codim}

\opn\Tr{Tr} \opn\bigrank{big\,rank}
\opn\superheight{superheight}\opn\lcm{lcm}
\opn\trdeg{tr\,deg}
\opn\reg{reg} \opn\lreg{lreg} \opn\ini{in} \opn\lpd{lpd}
\opn\size{size}\opn{\mult}{mult}
\opn\Div{Div} \opn\cl{cl} \opn\Cl{Cl}
\opn\Spec{Spec} \opn\Supp{Supp} \opn\supp{supp} \opn\Sing{Sing}
\opn\Ass{Ass} \opn\Min{Min}
\opn\Ann{Ann} \opn\Rad{Rad} \opn\Soc{Soc}
\opn\Syz{Syz} \opn\Im{Im} \opn\Ker{Ker} \opn\Coker{Coker}
\opn\Am{Am} \opn\Hom{Hom} \opn\Tor{Tor} \opn\Ext{Ext}
\opn\End{End} \opn\Aut{Aut} \opn\id{id}

\opn\nat{nat}
\opn\pff{pf}
\opn\Pf{Pf} \opn\GL{GL} \opn\SL{SL} \opn\mod{mod} \opn\ord{ord}
\opn\Gin{Gin}
\opn\Hilb{Hilb}\opn\adeg{adeg}\opn\std{std}\opn\ip{infpt}
\opn\Pol{Pol}
\opn\sat{sat}
\opn\Var{Var}
\opn\Gen{Gen}
\opn\vol{vol}

\opn\aff{aff} \opn\con{conv} \opn\relint{relint} \opn\st{st}
\opn\lk{lk} \opn\cn{cn} \opn\core{core} \opn\vol{vol}
\opn\link{link} \opn\star{star}
\opn\gr{gr}

\def\Cc{{\mathcal C}}

\def\Ic{{\mathcal I}}
\def\Gc{{\mathcal G}}
\def\Fc{{\mathcal F}}
\def\Pc{{\mathcal P}}
\def\Hc{{\mathcal H}}

\def\pot#1#2{#1[\kern-0.28ex[#2]\kern-0.28ex]}

\opn\dirlim{\underrightarrow{\lim}}
\opn\inivlim{\underleftarrow{\lim}}
\let\union=\cup

\def\Implies{\ifmmode\Longrightarrow \else
        \unskip${}\Longrightarrow{}$\ignorespaces\fi}
\def\implies{\ifmmode\Rightarrow \else
        \unskip${}\Rightarrow{}$\ignorespaces\fi}
\def\iff{\ifmmode\Longleftrightarrow \else
        \unskip${}\Longleftrightarrow{}$\ignorespaces\fi}

\let\:=\colon
\newtheorem{Theorem}{Theorem}[section]
\newtheorem{Lemma}[Theorem]{Lemma}
\newtheorem{Corollary}[Theorem]{Corollary}
\newtheorem{Proposition}[Theorem]{Proposition}

\newtheorem{Examples}[Theorem]{Examples}

\let\epsilon\varepsilon
\let\phi=\varphi
\let\kappa=\varkappa
\def\qed{\ifhmode\textqed\fi
      \ifmmode\ifinner\quad\qedsymbol\else\dispqed\fi\fi}
\def\textqed{\unskip\nobreak\penalty50
       \hskip2em\hbox{}\nobreak\hfil\qedsymbol
       \parfillskip=0pt \finalhyphendemerits=0}
\def\dispqed{\rlap{\qquad\qedsymbol}}

\opn\dis{dis}
\def\pnt{{\raise0.5mm\hbox{\large\bf.}}}

\opn\Lex{Lex}

\begin{document}

\title{
Shifted symmetric $\delta$-vectors of convex polytopes 
}

\author{Akihiro Higashitani}

\address{Akihiro Higashitani, 
Department of Pure and Applied Mathematics,
Graduate School of Information Science and Technology,
Osaka University,
Toyonaka, Osaka 560-0043, Japan}
\email{sm5037ha@ecs.cmc.osaka-u.ac.jp}

\begin{abstract}
Let $\Pc \subset \RR^N$ be an integral 
convex polytope of dimension $d$ and
$\delta(\Pc) 
= (\delta_0, \delta_1, \ldots, \delta_d)$ 
its Ehrhart $\delta$-vector.
It is known that
$
\sum_{j=0}^{i} \delta_{d-j} 
\leq
\sum_{j=0}^{i} \delta_{j+1}
$
for each $0 \leq i \leq [(d-1)/2]$.
A $\delta$-vector
$\delta(\Pc) 
= (\delta_0, \delta_1, \ldots, \delta_d)$
is called shifted symmetric
if 
$
\sum_{j=0}^{i} \delta_{d-j} 
=
\sum_{j=0}^{i} \delta_{j+1}
$
for each $0 \leq i \leq [(d-1)/2]$,
i.e.,
$\delta_{d-i} = \delta_{i+1}$
for each $0 \leq i \leq [(d-1)/2]$.
In this paper, some properties of integral convex polytopes 
with shifted symmetric $\delta$-vectors will be studied. 
Moreover, as a natural family of those, 
$(0,1)$-polytopes will be introduced. 
In addition, shifted symmetric $\delta$-vectors with (0,1)-vectors 
are classified when $\sum_{i=0}^d\delta_i \leq 5$. 
\end{abstract}

\maketitle
\section*{Introduction}
An \textit{integral} convex polytope is a convex polytope 
any of whose vertices has integer coordinates.
Let $\Pc \subset \RR^N$ be an integral convex 
polytope of dimension $d$ and
\[
i(\Pc,n) = |n\Pc \cap \ZZ^N|,
\, \, \, \, \, 
n = 1, 2, 3, \ldots.
\]
Here $n\Pc = \{ n\alpha : \alpha \in \Pc \}$
and $|X|$ is the cardinality of a finite set $X$.
The systematic study of $i(\Pc,n)$ originated
in the work of Ehrhart \cite{Ehrhart}, who established
the following fundamental properties:
\begin{enumerate}
\item[(0.1)]
$i(\Pc,n)$ is a polynomial in $n$ of degree $d$.
(Thus in particular $i(\Pc,n)$ can be defined 
for {\em every} integer $n$.)
\item[(0.2)]
$i(\Pc,0) = 1$.
\item[(0.3)]
(loi de r\'eciprocit\'e)
$( - 1 )^d i(\Pc, - n) 
= |n(\Pc - \partial \Pc) \cap \ZZ^N|$
for every integer $n > 0$.
\end{enumerate}
We say that $i(\Pc,n)$ is the {\em Ehrhart polynomial}
of $\Pc$. 
We refer the reader to
\cite[pp. 235--241]{StanleyEC}
and \cite[Part II]{HibiRedBook}
for the introduction 
to the theory of Ehrhart polynomials. 

We define the sequence 
$\delta_0, \delta_1, \delta_2, \ldots$ of integers
by the formula
\begin{eqnarray}
\label{delta}
(1 - \lambda)^{d + 1}
\left[ 1 + \sum_{n=1}^{\infty} i(\Pc,n) \lambda^n \right]
= \sum_{i=0}^{\infty} \delta_i \lambda^i.
\end{eqnarray}
It follows from
the basic fact (0.1) and (0.2) on $i(\Pc,n)$
together with a fundamental result on generating
function (\cite[Corollary 4.3.1]{StanleyEC})
that $\delta_i = 0$ for every $i > d$.
We say that the sequence 
\[
\delta(\Pc) 
= (\delta_0, \delta_1, \ldots, \delta_d)
\]
which appears in Eq.\,(\ref{delta}) is the
{\em $\delta$-vector} of $\Pc$.
Alternate names of $\delta$-vector are for example 
{\em Ehrhart $\delta$-vector}, {\em Ehrhart $h$-vector} 
or {\em $h^*$-vector}. 

Thus $\delta_0 = 1$ and 
$\delta_1 = |\Pc \cap \ZZ^N| - (d + 1)$.
Let $\partial \Pc$ denote 
the boundary of $\Pc$ and
\[
i^*(\Pc,n) = |n(\Pc - \partial \Pc) \cap \ZZ^N|,
\, \, \, \, \, 
n = 1, 2, 3, \ldots.
\]
By using (0.3) one has
\begin{eqnarray}
\label{deltadual}
\sum_{n=1}^{\infty} i^*(\Pc,n) \lambda^n 
= \frac
{\sum_{i=0}^{d} \delta_{d-i} \lambda^{i+1}}
{(1 - \lambda)^{d + 1}}.
\end{eqnarray}
In particular, 
\[
\delta_d = |(\Pc - \partial \Pc) \cap \ZZ^N|.
\]
Hence $\delta_1 \geq \delta_d$. 
Moreover, each $\delta_i$ is nonnegative 
(\cite{StanleyDRCP}).
In addition, 
if $(\Pc - \partial \Pc) \cap \ZZ^N$ is nonempty,
then one has 
$\delta_1 \leq \delta_i$ for every 
$1 \leq i \leq d - 1$
(\cite{HibiLBT}).

When $d = N$, the leading coefficient 
$(\sum_{i=0}^{d}\delta_i)/d!$
of $i(\Pc,n)$ is equal to the usual volume of $\Pc$
(\cite[Proposition 4.6.30]{StanleyEC}).
In general, the positive integer 
$\vol(\Pc) = \sum_{i=0}^{d}\delta_i$
is said to be the {\em normalized volume} of $\Pc$.

It follows from Eq.\,(\ref{deltadual}) that 
\[
\max\{ i : \delta_i \neq 0 \}
= d+1- \min\{ i : 
i(\Pc - \partial \Pc) \cap \ZZ^N \neq \emptyset \}.
\] 

Recently, $\delta$-vectors of integral convex polytopes 
have been studied intensively. 
(For example, see \cite{Payne},\cite{Staple1} and \cite{Staple2}.) 

There are two well-known inequalities of $\delta$-vectors. 
Let $s = \max\{ i : \delta_i \neq 0 \}$.
Stanley \cite{StanleyJPAA} shows the inequalities 
\begin{eqnarray}
\label{Stanley}
\delta_0 + \delta_1 + \cdots + \delta_i
\leq \delta_s + \delta_{s-1} + \cdots + \delta_{s-i},
\, \, \, \, \, 
0 \leq i \leq [s/2]
\end{eqnarray}
by using the theory of Cohen--Macaulay rings.
On the other hand, the inequalities  
\begin{eqnarray}
\label{Hibi}
\delta_{d} + \delta_{d-1} + \cdots + \delta_{d-i}
\leq \delta_1 + \delta_2 + \cdots + \delta_i
+ \delta_{i+1},
\, \, \, \, \, 
0 \leq i \leq [(d-1)/2]
\end{eqnarray}
appear in \cite[Remark (1.4)]{HibiLBT}.
The above two inequalities are generalized in \cite{Staple1}. 

A $\delta$-vector 
$\delta(\Pc) 
= (\delta_0, \delta_1, \ldots, \delta_d)$
is called {\em symmetric}
if the equalities hold 
in Eq. (\ref{Stanley}) for each 
$0 \leq i \leq [s/2]$, i.e.,
$\delta_i = \delta_{s-i}$
for each $0 \leq i \leq [s/2]$.
The $\delta$-vector $\delta(\Pc)$
of $\Pc$ is symmetric if and only if
the Ehrhart ring \cite[Chapter X]{HibiRedBook}
of $\Pc$ is Gorenstein.  
A combinatorial characterization
for the $\delta$-vector to be symmetric
is studied in \cite{HibiCombinatorica}
and \cite{DeNegriHibi}.

We say that a $\delta$-vector 
$\delta(\Pc) 
= (\delta_0, \delta_1, \ldots, \delta_d)$
is {\em shifted symmetric} if the equalities hold
in Eq. (\ref{Hibi}) for each 
$0 \leq i \leq [(d-1)/2]$,
i.e.,
$\delta_{d-i} = \delta_{i+1}$
for each $0 \leq i \leq [(d-1)/2]$.
It seems likely that an integral convex polytope
with a shifted symmetric $\delta$-vector is quite rare.
Thus it is reasonable to sutudy a property of and 
to find a natural family of 
integral convex polytopes with shifted symmetric $\delta$-vectors. 
In section 2, some characterizations of an integral convex polytope 
with a shifted symmetric $\delta$-vector are given. 
Concretely, in Theorem \ref{faces1}, it is shown that 
integral convex polytopes with shifted symmetric $\delta$-vectors 
have a special property. 
Moreover, 
as a generalization of an integral convex polytope 
with a shifted symmetric $\delta$-vector, 
an integral simplicial polytope 
any of whose facet has the normalized volume 1 
is considered in section 2. 
In section 3, 
a family of $(0,1)$-polytopes 
with shifted symmetric $\delta$-vectors is presented. 
These shifted symmetric $\delta$-vectors are $(0,1)$-vectors. 
In addition, by using those examples, 
we classify shifted symmetric $\delta$-vectors with (0,1)-vectors 
when $\sum_{i=0}^d\delta_i \leq 5$ in section 4.

\section{Review on the computation of the $\delta$-vector of a simplex}

We recall from \cite[Part II]{HibiRedBook}
the well-known combinatorial technique how to compute
the $\delta$-vector of a simplex. 

\begin{itemize}
\item
Given an integral $d$-simplex $\Fc \subset \RR ^N$ with the vertices
$v_0, v_1, \ldots, v_d$, we set 
$\widetilde \Fc=\left\{(\a,1)\in \RR^{N+1} \, : \, \a \in \Fc \right\}$. 
And
$
\partial \widetilde \Fc=\left\{(\a,1)\in \RR^{N+1} \, : \, \a \in \partial \Fc \right\}
$
is its boundary.
\item
Let $\Cc(\widetilde \Fc)=\Cc=
\{r \beta \, : \, \beta \in \widetilde \Fc,0 \leq r \in \QQ \} 
\subset \RR^{N+1}.$ 
Its boundary is
$
\partial \Cc = \left\{r\beta \, : \,  \beta \in \partial \widetilde \Fc,0 \leq r \in \QQ \right\}.
$
\item
Let $S$ (resp. $S^*$) be the set of all points $\a \in \Cc \cap \ZZ^{N+1}$
(resp. $\a \in (\Cc- \partial \Cc) \cap \ZZ^{N+1}$)
of the form
$
\a= \sum_{i=0}^{d}r_i(v_i,1),
$
where each $r_i \in \QQ$ with $0 \leq r_i<1$ (resp. with $0<r_i \leq 1$).
\item
The degree of an integer point $(\a,n) \in \Cc$ is
$
\deg(\a,n):=n.
$
\end{itemize}

\begin{Lemma}
\label{computation}
{\em (a)}
Let $\delta_i$ be the number of integer points $\a \in S$
with $\deg \a=i$.  Then
\[
1+ \sum_{n=1}^{\infty}i(\Fc,n) \lambda^n=\frac{\delta_0+\cdots+\delta_d\lambda^d}{(1-\lambda)^{d+1}}.
\]

{\em (b)}
Let $\delta_i^*$ be the number of integer points $\a \in S^*$
with $\deg\a=i$.  Then
\[
\sum_{n=1}^{\infty}i^*(\Fc,n) \lambda^n=\frac{\delta_1^*\lambda+\cdots+\delta_{d+1}^*\lambda^{d+1}}{(1-\lambda)^{d+1}}.
\]

{\em (c)}
One has $\delta_i^* = \delta_{(d+1)-i}$ for each $1 \leq i \leq d+1$.
\end{Lemma}

We say that a $\delta$-vector 
$\delta(\Pc) 
= (\delta_0, \delta_1, \ldots, \delta_d)$
is {\em shifted symmetric} 
if $\delta_{d-i} = \delta_{i+1}$ 
for each $0 \leq i \leq [(d-1)/2]$. 
Since $\delta_d = \delta_1$, 
an integral convex polytope with a shifted symmetric $\delta$-vector 
is always a $d$-simplex. 

The followings are some examples of a simplex with a shifted symmetric $\delta$-vector. 

Let ${\bf e}_i$ denote
the $i$th
canonical unit coordinate vector of $\RR^d$. 

\begin{Examples}\label{example}
{\em (a) 
We define $v_i \in \RR^d,i=0,1,\ldots,d$ by setting 
$v_i = {\bf e}_i$ for $i=1,\ldots,d$ and $v_0 = (-e,\ldots,-e)$, 
where $e$ is a nonnegative integer. 
Let $\Pc=\con\{v_0,v_1,\ldots,v_d\}$. 
Then one has $\vol(\Pc) = ed + 1$ 
by using an elementary linear algebra. 
When $e=0$, it is clear that $\delta(\Pc)=(1,0,0,\ldots,0)$. 
When $e$ is positive, we know that 
$$
\frac{j}{ed+1}\sum_{i=1}^d(v_i,1)+\frac{(e-j)d+1}{ed+1}(v_0,1) = (j-e,j-e,\ldots,j-e,1) 
$$
and $0 < \frac{j}{ed+1}, \frac{(e-j)d+1}{ed+1} < 1$ for every $1 \leq j \leq e$. 
Then Lemma \ref{computation} says that $\delta_1,\delta_d \geq e$. 
Since $\delta_i \geq \delta_1$ for $1 \leq i \leq d-1$ and $\vol(\Pc)=ed+1$, 
we obtain $\delta(\Pc)=(1,e,e,\ldots,e)$. 

(b) 
Let $d \geq 3$. We define $v_i \in \RR^d,i=0,1,\ldots,d$ by setting 
$v_i = {\bf e}_i$ for $i=1,\ldots,d$ and $v_0 = (e,\ldots,e)$, 
where $e$ is a positive integer. 
Let $\Pc=\con\{v_0,v_1,\ldots,v_d\}$. 
Then one has $\vol(\Pc) = ed - 1$ 
by using an elementary linear algebra. 
And we know that 
$$
\frac{j}{ed-1}\sum_{i=1}^d(v_i,1)+\frac{(e-j)d-1}{ed-1}(v_0,1) = (e-j,e-j,\ldots,e-j,1) 
$$
and $0 < \frac{j}{ed-1}, \frac{(e-j)d-1}{ed-1} < 1$ 
for every $1 \leq j \leq e-1$. 
Thus $\delta_1,\delta_d \geq e-1$ by Lemma \ref{computation}. 
In addition we know that 
$$
\frac{ke+j}{ed-1}\sum_{i=1}^d(v_i,1)+\frac{(e-j)d-1-k}{ed-1}(v_0,1) 
= (e-j,e-j,\ldots,e-j,k+1) 
$$
and $0 <\frac{ke+j}{ed-1}, \frac{(e-j)d-1-k}{ed-1}< 1$ 
for every $0 \leq j \leq e-1$ and $1 \leq k  \leq d-2$. 
Hence $\delta(\Pc)=(1,e-1,e,e,\ldots,e,e-1)$. 
}
\end{Examples}

\section{Some characterizations of an integral convex polytope 
with a shifted symmetric $\delta$-vector}

In the first half of this section, 
two results of an integral convex polytope 
with a shifted symmetric $\delta$-vector are given. 
And in the latter half of this section, 
we generalize an integral convex polytope with a shifted symmetric $\delta$-vector. 

\begin{Theorem}\label{faces1}
Let $\Pc$ be a $d$-simplex whose vertices are 
$v_0,v_1,\ldots,v_d \in \RR^d$ 
and $S \; (S^*)$ the set which appears in section $1$. 
Then the following conditions are equivalent: 
\begin{enumerate}
\item[(i)]
$\delta(\Pc)$ is shifted symmetric; 
\item[(ii)]
the normalized volume of all facets of $\Pc$ is equal to $1$; 
\item[(iii)]
each element $(\a,n) \in S \backslash \{(0,\ldots,0,0)\}$ 
has a unique expression of the form: 
\begin{equation}\label{form}
(\a,n)=\sum_{j=0}^{d}r_j(v_j,1) \;\; \text{with} \;\; 0<r_j<1 \;\; \text{for} \;\; j=0,1,\ldots,d, 
\end{equation}
where $\a \in \ZZ^d$ and $n \in \ZZ$. 
\end{enumerate}
\end{Theorem}

\begin{proof}
{\bf ((i) $\Leftrightarrow$ (iii))}
If each element $(\a,n) \in S \backslash \{(0,\ldots,0,0)\}$ 
has the form (\ref{form}), 
each element 
$(\alpha',n') \in S^* \backslash \{(\sum_{j=0}^dv_j,\ldots,\sum_{j=0}^dv_j,d+1)\}$ 
also has the same form (\ref{form}). 
This implies that $\d(\Pc)$ is shifted symmetric. 
On the other hand, suppose that $\delta(\Pc)$ is shifted symmetric. 
Let $\min\{i:\delta_i \not= 0, i>0 \} = s_1$ and $\delta_{s_1}=m_1(\not=0)$. 
Then one has $d+1-\max\{i:\delta_i \not= 0 \} = s_1$ and 
both $S$ and $S^*$ have the $m_1$ elements with degree $s_1$. 
If an element 
$(\alpha',s_1) \in S^*$ does not have the form (\ref{form}), 
there is $0 \leq j \leq d$ with $r_j=1$, 
say, $r_0=r_1=\cdots=r_a=1$ and $0 < r_{a+1},r_{a+2},\ldots,r_d < 1$. 
Then $S$ has an element 
$(\a'-v_0-v_1-\cdots-v_a,s_1-a-1)\not=(0,\ldots,0,0)$, a contradiction. 
Thus each element $(\a',s_1) \in S^*$ has the form (\ref{form}). 
If we set $\min\{i:\delta_i\not=0,i>s_1 \} = s_2$, 
the same discussions can be done as above. 
Thus each element $(\b',s_2) \in S^*$ has the form (\ref{form}). 
Hence each element $(\alpha',n') \in S^* \backslash \{(\sum_{j=0}^dv_j,\ldots,\sum_{j=0}^dv_j,d+1)\}$ 
has the form (\ref{form}), that is to say, 
each element $(\a,n) \in S \backslash \{(0,\ldots,0,0)\}$ has the form (\ref{form}).

{\bf ((ii) $\Leftrightarrow$ (iii))}
Let $\delta(\Pc)=(\delta_0,\delta_1,\ldots,\delta_d) \in \ZZ^{d+1}$ 
be the $\delta$-vector of $\Pc$ and 
$\delta(\Fc)=(\delta_0',\delta_1',\ldots,\delta_{d-1}') \in \ZZ^d$ 
the $\delta$-vector of a facet $\Fc$ of $\Pc$. 
Then one has $\delta_i' \leq \delta_i$ for $0 \leq i \leq d-1$. 
If there is a facet $\Fc$ with $\vol(\Fc)\not=1$, 
say, its vertices are $v_0,v_1,\ldots,v_{d-1}$, 
there exists an element $(\a,n) \in S$ 
with $\a=\sum_{j=0}^{d-1}r_jv_j+0 \cdot v_d$ and $n > 0$. 
This implies that there exists an element of $S \backslash \{(0,\ldots,0,0)\}$ 
which does not have the form (\ref{form}). 
On the other hand, suppose that 
there exists an element $(\a,n) \in S \backslash \{(0,\ldots,0,0)\}$ 
which does not have the form (\ref{form}), 
i.e., $(\a,n)=\sum_{j=0}^{d}r_j(v_j,1)$ 
and there is $0 \leq j \leq d$ with $r_j=0$, say, $r_d=0$. 
Then the normalized volume of the facet whose vertices are $v_0,v_1,\ldots,v_{d-1}$ 
is not equal to 1. 
\end{proof}




\begin{Theorem}
Let $\Pc$ be a $d$-simplex. 
If $\vol(\Pc)=p$ with a prime number $p$ and 
$$
\min\{ i:\delta_i \not=0 ,i>0 \} = d+1- \max\{ i:\delta_i \not= 0 \}, 
$$
then $\delta(\Pc)$ is shifted symmetric. 
\end{Theorem}

\begin{proof}

The elements of $S$ form a cyclic group of prime order. 
Then every non-zero element of $S$ is a generator. 
Thus it can be considered 
whether $S \backslash \{(0,\ldots,0,0)\}$ is disjoint from $\partial S$, 
which case satisfies the condition of Theorem \ref{faces1} (iii), 
or $S$ is contained in a facet of $\partial S$, 
where $\partial S$ is a cyclic group generated by the vertices of a facet of $\Pc$. 
In the latter case, let $x \in S$ be an element of the maximal degree deg$(x)$, 
and let $-x$ denote its inverse. 
Since $S$ is contained in $\partial S$, deg$(x) + $deg$(-x) < d + 1$, 
which contradicts the assumption. 
Therefore, since Theorem \ref{faces1}, $\d(\Pc)$ is shifted symmetric. 
\end{proof}

Recall that an integral convex polytope with a shifted symmetric $\delta$-vector 
is always a simplex. 
Then we expand the definition of shifted symmetric to an integral simplicial polytope. 
We study an integral simplicial polytope $\Pc$ 
any of whose facet has the normalized volume 1. 
When $\Pc$ is a simplex, its $\delta$-vector is shifted symmetric by Theorem \ref{faces1}. 

Let $h(\Delta(\Pc))=(h_0,h_1,\ldots,h_d)$ 
denote the $h$-vector of the boundary complex of $\Pc$. 
(See, \cite[Part I]{HibiRedBook}.) 
Then the following is a well-known fact 
about a lower bound of the $h$-vector for a simplicial $(d-1)$-sphere. 

\begin{Lemma}\label{lb}{\em(\cite{bar1},\cite{bar2})}
The $h$-vector $h(\Delta(\Pc))=(h_0,h_1,\ldots,h_d)$ 
of a simplicial $(d-1)$-sphere satisfies 
$h_1 \leq h_i$ for every $1 \leq i \leq d-1$. 
\end{Lemma}

Now, all of $h$-vectors of simplicial $(d-1)$-spheres 
satisfying the lower bound, i.e., $h_1 = h_i$ for every $1 \leq i \leq d-1$, 
are given by 
$h$-vectors of the boundary complexes of simplicial polytopes 
any of whose facet has the normalized volume 1. 
In fact, 
\begin{Theorem}
For an arbitrary positive integer $h_1$, 
there exists a $d$-dimensional simplicial polytope $\Pc$ 
any of whose facet has the normalized volume 1 and 
whose $h$-vector of the boundary complex coincides with 
$(1,h_1,\ldots,h_1,1) \in \ZZ^{d+1}$. 
\end{Theorem}

\begin{proof}
Let $d=2$. A convex polygon is always simplicial. 
And each facet of an integral polygon has the normalized volume 1 
if and only if there is no integer point in its boundary except its vertices. 
Hence, for an arbitrary positive integer $h_1$, 
we can say that there exists an integral polygon with $h_1 + 2$ vertices 
any of whose facet has the normalized volume 1. 

We assume when $d \geq 3$. 
Let $\Pc$ be the $d$-dimensional integral convex polytope 
whose vertices $v_i \in \RR^d$, $i=0,1,\ldots,d+n$, are of the form: 
\begin{eqnarray*}
v_i=
\begin{cases}
(0,\ldots,0) \;\;\;\; &\text{for} \;\;\;i=0, \\
{\bf e}_i &\text{for} \;\;\;i=1,2,\ldots,d, \\
(c_j,\ldots,c_j,j) &\text{for} \;\;\;i=d+1,d+2,\ldots,d+n, 
\end{cases}
\end{eqnarray*}
where $n$ is a positive even number, $j=i-d$ and $c_j=n+\frac{(n-j)(j-1)}{2}$. 

\textit{First step.} 
We prove that $\Pc$ is a simplicial convex polytope. 
We define the $\{(n+1)(d-1)+2\}$ convex hulls by setting 
$$
\begin{cases}
\Fc_i:=\con\{v_0,v_1,\ldots,v_{i-1},v_{i+1},\ldots,v_{d}\} \;\;\;\; \text{for} \;\;\;i=1,\ldots,d, \\
\Fc':=\con\{v_1,\ldots,v_{d-1},v_{d+1}\}, \\
\Gc_{i,j}:=\con\{v_1,\ldots,v_{i-1},v_{i+1},\ldots,v_{d-1},v_{d+j},v_{d+j-1}\} \;\;\;\;\; 
\text{for} \;\;\; i=1,\ldots,d-1,j=2,\ldots,n, \\
\Gc'_i:=\con\{v_1,\ldots,v_{i-1},v_{i+1},\ldots,v_{d-1},v_d,v_{d+n}\} \;\;\;\;\; \text{for} \;\;\; i=1,\ldots,d-1, 
\end{cases}
$$
and the followings are the equations of the hyperplanes containing the above convex hulls: 
$$
\begin{cases}
\Hc_i \supset \Fc_i:  -x_i=0 \;\;\; \text{for} \;\;\; i=1,\ldots,d, \\
\Hc' \supset \Fc': \sum_{k=1}^{d-1}x_k-(n(d-1)-1)x_d=1, \\
\Ic_{i,j} \supset \Gc_{i,j}: c'_j\sum_{k=1}^{i-1}x_k-(1-(d-2)c'_j)x_i+
c'_j\sum_{k=i+1}^{d-1}x_k+(j-\frac{n+2}{2})x_d =c_j' \\ 
\text{for} \;\;\; i=1,\ldots,d-1,j=2,\ldots,n, \;\;\; \text{where} \;\;\; 
c'_j=((j-1)c_j-jc_{j-1})=\frac{j^2-j+n}{2}, \\
\Ic'_i \supset \Gc'_i: n\sum_{k=1}^{i-1}x_k-(n(d-1)-1)x_i+n\sum_{k=i+1}^dx_k=n 
\;\;\; \text{for} \;\;\; i=1,\ldots,d-1. 
\end{cases}
$$

We prove that these $\{(n+1)(d-1)+2\}$ convex hulls are all facets of $\Pc$. 
If we write $\Hc \subset \RR^d$ for the hyperplane 
defined by the equation $a_1x_1+\cdots+a_dx_d = b$, 
then we write $\Hc^{(+)} \subset \RR^d$ for the closed half-space 
defined by the inequality $a_1x_1+\cdots+a_dx_d \leq b$. 

\begin{itemize}
\item
Let $\Pc_{n+1}=\con\{ v_0,\ldots,v_d \}$. 
Then one has $\Pc_{n+1} = (\bigcap_{i=1}^d\Hc_i^{(+)}) \cap (x_1+\cdots+x_d \leq 1)$. 
\item
Let $\Pc_k=\con\{ \Pc_{k+1} \union \{v_{d+k}\} \}$ for $k=n,n-1,\ldots,1$. 
Then it can be shown easily that 
$$\Pc_k = (\bigcap_{i=1}^d\Hc_i^{(+)}) \cap (\bigcap_{i=1}^{d-1} \Ic_i^{(+)'}) 
\cap (\bigcap_{\begin{subarray}{c}1 \leq i \leq d-1 \\ k+1 \leq j \leq n \end{subarray}}
\Ic_{i,j}^{(+)}) \cap 
(kx_1+\cdots+kx_{d-1}-(c_k(d-1)-1)x_d \leq k). $$
\end{itemize}
Then one has $\Pc_1=\con\{ \Pc_2 \union \{v_{d+1}\} \}=\Pc$. 
Thus we obtain the following equality: 
$$
\Pc=(\bigcap_{i=1}^d\Hc_i^{(+)}) \cap (\bigcap_{i=1}^{d-1}\Ic_i^{(+)'}) \cap 
(\bigcap_{\begin{subarray}{c}1 \leq i \leq d-1 \\ 2 \leq j \leq n \end{subarray}} 
\Ic_{i,j}^{(+)})\cap \Hc'^{(+)}. 
$$
Hence we can say that 
$\Fc_i,\Fc',\Gc_{i,j}$ and $\Gc_i'$ are all facets of $\Pc$ 
and they are $(d-1)$-simplices. 

\textit{Second step.} 
We prove that the normalized volume of each facet of $\Pc$ 
is equal to 1. 
One has $\vol(\Fc_i)=1$ for $i=1,\ldots,d$ since $\vol(\Pc_{n+1})=1$ 
and one has $\vol(\Gc_i')=1$ for $1 \leq i \leq d-1$ 
since $\con\{v_1,\ldots,v_d,v_{d+n}\}$ is a simplex 
with a shifted symmetric $\delta$-vector by Examples \ref{example}(b). 
And one has $\vol(\Fc')=1$ 
since $\vol(\con\{v_0,v_1,\ldots,v_{d-1},v_{d+1}\})=1$. 
When we consider $\vol(\Gc_{i,j})$, 
we are enough to prove that $\vol(\Gc_{d-1,j})=1$ by the symmetry. 

For a $(d-1)$-simplex $\Gc_{d-1,j}$, 
we consider the elements of the set $S$ which appears in section $1$: 
$$
(\a_1,\ldots,\a_d,r)=
r_1(v_1,1)+\cdots+r_{d-2}(v_{d-2},1)+r_{d-1}(v_{d+j},1)+r_d(v_{d+j-1},1), 
$$
where $(\a_1,\ldots,\a_d,r) \in \ZZ^{d+1}$ and $0 \leq r_i < 1$ for $0 \leq i \leq d$. 
Then one has 
$$
(\a_1,\ldots,\a_d) = 
(r_1+r_{d-1}c_j+r_dc_{j-1},\ldots,r_{d-2}+r_{d-1}c_j+r_dc_{j-1},
r_{d-1}c_j+r_dc_{j-1},r_{d-1}j+r_d(j-1)). 
$$
Since $r_1=\a_1-\a_{d-1} \in \ZZ$, we obtain $r_1=0$. 
Similary we obtain $r_2=\cdots=r_{d-2}=0$. 
Hence we can rewrite 
$(\a_1,\ldots,\a_d,r)=r_{d-1}(v_{d+j},1)+r_d(v_{d+j-1},1)$. 
It then follows that 
$$
\vol(\Gc_{d-1,j})=\vol(\con\{v_{d+j},v_{d+j-1}\})=
\vol(\con\{(c_j,j),(c_{j-1},j-1)\}) 
$$
Since $r_{d-1}j+r_d(j-1) \in \ZZ$ and $r_{d-1}+r_d \in \ZZ$, one has $r_{d-1}=r_d=0$. 
This implies that $\vol(\con\{(c_j,j),(c_{j-1},j-1)\})=1$. 
Thus $\vol(\Gc_{d-1,j})=1$. 

\textit{Third step.} 
By the first step and the second step, 
$\Pc$ is a $d$-dimensional simplicial polytope 
with $\{(n+1)(d-1)+2\}$ facets and $(d+n+1)$ vertices 
any of whose facet has the normalized volume 1. 
Hence, by Lemma \ref{lb}, one has 
$h(\Delta(\Pc))=(1,n+1,\ldots,n+1,1)$ for a positive even number $n$. 
Thus, when $h_1$ is odd and $h_1 \geq 3$, we know that 
there exists a simplicial polytope with $h(\Delta(\Pc))=(1,h_1,\ldots,h_1,1)$ 
any of whose facet has the normalized volume 1. 
When $h_1=1$, it is clear that $h(\Delta(\Pc_{n+1}))=(1,1,\ldots,1)$. 
When $h_1$ is even and $h_1 \geq 2$, let $\Pc'=\con\{v_1,\ldots,v_{d+n}\}$. 
Then we can verify that $\Pc'$ is 
a simplicial polytope with $h(\Delta(\Pc'))=(1,n,\ldots,n,1)$ 
any of whose facet has the normalized volume 1. 
\end{proof}

\section{A family of $(0,1)$-polytopes with shifted symmetric $\delta$-vectors}

In this section, 
a family of $(0,1)$-polytopes with shifted symmetric $\delta$-vectors 
is studied. 
We classify completely the $\delta$-vectors of those polytopes. 
Moreover, we consider when those $\delta$-vectors 
are both symmetric and shifted symmetric. 

Let $d=m+n$ with positive integers $m$ and $n$. 
We study the $\delta$-vector 
of the integral convex polytope $\Pc \subset \RR^d$ 
whose vertices are of the form: 
\begin{eqnarray}
\label{vertex}
v_i = 
\begin{cases}
{\bf e}_i + {\bf e}_{i+1} + \cdots + {\bf e}_{i+m-1} \;\;\;\; &i=1,\ldots,d, \\
(0,\ldots,0) &i=0, 
\end{cases}
\end{eqnarray}
where ${\bf e}_{d+i} = {\bf e}_i$. 

The normalized volume of $\Pc$ is equal to 
the absolute value of the determinant of the circulant matrix 
\begin{eqnarray}
\label{matrix}
\begin{vmatrix}
v_1 \\
\vdots \\
v_{d} 
\end{vmatrix}. 
\end{eqnarray}
This determinant (\ref{matrix}) can be calculated easily. 
In fact, 
\begin{Proposition}\label{m}
When $(m,n)=1$, the determinant (\ref{matrix}) is equal to $\pm m$. 
And when $(m,n)\not=1$, the determinant (\ref{matrix}) is equal to $0$. 
Here $(m,n)$ is the greatest common divisor of $m$ and $n$. 
\end{Proposition}
A proof of this proposition can be given 
by the formula of the determinant of the circulant matrix. 
Thus one has $\vol(\Pc)=m$ when $(m,n)=1$. 

In this section, we assume only the case of $(m,n)=1$. 
Hence $(m,d)=1$. 

For $j = 1,2,\ldots,d-1,$ 
let $q_j$ be the quotient of $jm$ divided by $d$ 
and $r_j$ its remainder 
i.e., one has the equalities
$$
jm=q_jd + r_j \quad \text{for} \quad j=1,2,\ldots,d-1. 
$$
It then follows from $(m,d)=1$ that 
$$
0 \leq q_j \leq m-1, 1 \leq r_j \leq d-1 
$$
and 
$$
r_j \not= r_{j'} \text{ if } j\not=j' 
$$
for every $1 \leq j,j' \leq d-1$. 
In addition, for $k=1,2,\ldots,m-1$, 
let $j_k \in \{1,2,\ldots.d-1\}$ be the integer with $r_{j_k}=k$, 
i.e., one has the equalities
$$
j_km=q_{j_k}d+r_{j_k}=q_{j_k}d+k \quad \text{for} \quad k=1,2,\ldots,m-1. 
$$
Then $q_{j_k} > 0$. 
Thus one has 
$$
1 \leq q_{j_k},r_{j_k} \leq m-1 
$$
for every $1 \leq k \leq m-1$. 

For an integer $a$, let $\overline{a}$ denote the residue class in $\ZZ/d\ZZ$. 

\begin{Theorem}\label{main}
Let $\Pc$ be the integral convex polytope 
whose vertices are of the form (\ref{vertex}) 
and $\d(\Pc)=(\delta_0,\delta_1,\ldots,\delta_d)$ its $\delta$-vector. 
For each $1 \leq i \leq d$, 
one has $\overline{im} \in \{ \overline{1},\overline{2},\ldots,\overline{m-1} \}$ 
if and only if one has $\delta_i=1$. 
Moreover, $\d(\Pc)$ is shifted symmetric, 
i.e., $\delta_{i+1}=\delta_{d-i}$ for each $0 \leq i \leq [(d-1)/2]$. 
\end{Theorem}

\begin{proof}
By using the above notations, 
we obtain 
\begin{align*}
\frac{q_{j_k}}{m}\left\{(v_1,1)+(v_2,1)+\cdots+(v_d,1)\right\}+\frac{r_{j_k}}{m}(v_0,1) 
&=(q_{j_k},\ldots,q_{j_k},j_k) \in \ZZ^{d+1} 
\end{align*} 
and $0 < \frac{q_{j_k}}{m},\frac{r_{j_k}}{m} < 1$ 
for every $1 \leq k \leq m-1$. 
Then Lemma \ref{computation} guarantees that 
one has $\delta_{j_k} \geq 1$ for $k=1,\ldots,m-1$. 
Considering $\sum_{i=0}^d\delta_i=m$ by Proposition \ref{m}, it turns out that 
$\delta(\Pc)$ coincides with 
\begin{eqnarray*}
\delta_i=
\begin{cases}
1 \qquad i=0,j_1,j_2,\ldots,j_{m-1}, \\
0 \qquad otherwise. 
\end{cases}
\end{eqnarray*}
Now $\overline{im} \in \{ \overline{1},\overline{2},\ldots,\overline{m-1} \}$ 
is equivalent 
with $i \in \{j_1,\ldots,j_{m-1}\}$. 
Therefore one has $\delta_i=1$ if and only if 
$\overline{im} \in \{ \overline{1},\overline{2},\ldots,\overline{m-1} \}$ 
for each $1 \leq i \leq d$. 

In addition, by virtue of Theorem \ref{faces1}, 
$\d(\Pc)$ is shifted symmetric, as required. 

\end{proof}

\begin{Corollary}\label{corollary}
Let $\Pc$ be the integral convex polytope 
whose vertices are of the form (\ref{vertex}) 
and $\d(\Pc)=(\delta_0,\delta_1,\ldots,\delta_d)$ its $\delta$-vector. 
Then $\d(\Pc)$ is symmetric, 
i.e., $\delta_i=\delta_{s-i}$ for each $0 \leq i \leq [s/2]$ 
if and only if one has $d \equiv m-1 \pmod{m}$. 
\end{Corollary}
\begin{proof}
Let $p$ be the quotient of $d$ divided by $m$ and $r$ its remainder, 
i.e., one has $d=mp+r$. 
And let $j_t=\min\{j_1,j_2,\ldots,j_{m-1}\}$. 
On the one hand, one has $j_tm=d+t$. 
On the other hand, one has $(p+1)m=d+m-r$ and $1 \leq m-r \leq m-1$. 
It then follows from Theorem \ref{main} that 
$p+1=j_t=\min\{i:\delta_i\not=0,i>0\}$. 
Hence $d-p = s = \max\{ i : \delta_i \neq 0 \}$ 
since $\d(\Pc)$ is shifted symmetric. 

When $d \equiv m-1 \pmod{m}$, i.e., $r=m-1$, 
we can obtain the equalities 
\begin{eqnarray*}
d-p=mp+r-p&=&mp+m-1-p=(m-1)(p+1). 
\end{eqnarray*}
In addition, for nonnegative integers $l(p+1)$, $l=1,2,\ldots,m-1$, 
the following equalities hold: 
$$
\overline{l(p+1)m}=\overline{l(mp+m)}=\overline{l(mp+m-1)+l}=\overline{ld+l}=
\overline{l} \in \{ \overline{1},\overline{2},\ldots,\overline{m-1} \}. 
$$
Thus it turns out that $\d(\Pc)$ coincides with 
\begin{eqnarray*}
\delta_i=
\begin{cases}
1 \quad i=0,p+1,2(p+1),\ldots,(m-1)(p+1), \\
0 \quad otherwise, 
\end{cases}
\end{eqnarray*}
by Theorem \ref{main}. 
It then follows that 
$$
\delta_{k(p+1)}=\delta_{(m-1-k)(p+1)}=\delta_{s-k(p+1)}=1 
$$ 
for every $0 \leq k \leq m-1$ and 
$$
\delta_i=\delta_{s-i}=0
$$ 
for every $0 \leq i \leq s$ with $i \not= k(p+1)$, $k=0,1,\ldots,m-1$. 
These equalities imply that $\d(\Pc)$ is symmetric. 

Suppose that $\d(\Pc)$ is symmetric. 
Our work is to show that $r=m-1$. 
Then one has 
$$
\delta_0=\delta_s=\delta_{d-p}=\delta_{(m-1)(p+1)}=1. 
$$
Since $\d(\Pc)$ is also shifted symmetric, 
one has $\delta_{(m-1)(p+1)}=\delta_{p+1}$. 
Hence one has 
$\delta_{p+1}=\delta_{(m-2)(p+1)}=\delta_{2(p+1)}=\cdots=\delta_{[(m-1)/2](p+1)}$=1 
since $\d(\Pc)$ is both symmetric and shifted symmetric. 
When $m$ is odd, one has $\frac{d-p}{2}=\frac{m-1}{2}(p+1)$ 
since $\d(\Pc)$ is symmetric. Thus $r=m-1$. 
When $m$ is even, one has $\frac{d+1}{2}=\frac{m}{2}(p+1)$ 
since $\d(\Pc)$ is shifted symmetric. Thus $r=m-1$. 

Therefore $\d(\Pc)$ is symmetric 
if and only if $d \equiv m-1 \pmod{m}$. 
\end{proof}

\section{Classifications of shifted symmetric $\d$-vectors 
with (0,1)-vectors}

In this section, we will classify 
all the possible shifted symmetric symmetric $\d$-vectors 
with (0,1)-vectors when $\sum_{i=0}^d\delta_i \leq 5$ 
by using the examples in the previous section. 

In \cite{smallvolume}, the possible $\delta$-vectors of 
integral convex polytopes are classified completely 
when $\sum_{i=0}^d\delta_i \leq 3$. 

\begin{Lemma}\label{small}
Let $d \geq 3$.
Given a finite sequence
$(\delta_0, \delta_1, \ldots, \delta_d)$
of nonnegative integers, where $\delta_0 = 1$
and $\delta_1 \geq \delta_d$,
which satisfies $\sum_{i=0}^{d} \delta_i \leq 3$,
there exists an integral convex
polytope $\Pc \subset \RR^d$
of dimension $d$
whose $\delta$-vector
coincides with
$(\delta_0, \delta_1, \ldots, \delta_d)$
if and only if
$(\delta_0, \delta_1, \ldots, \delta_d)$
satisfies all inequalities
{\em (\ref{Stanley})} and {\em (\ref{Hibi})}.
\end{Lemma}

As an analogy of Lemma \ref{small}, 
we classify shifted symmetric symmetric $\delta$-vectors 
with (0,1)-vectos when $\sum_{i=0}^d\delta_i =4$ or 5. 

Now, in what follows, a sequence
$(\delta_0, \delta_1, \ldots, \delta_d)$
with each $\delta_i \in \{ 0, 1 \}$,
where $\delta_0 = 1$, 
which satisfies all inequalities
(\ref{Stanley}) and 
all equalities (\ref{Hibi})
together with
$\sum_{i=0}^{d} \delta_i = 4$ or 5 
will be considered.

At first, we consider the case of $\sum_{i=0}^{d} \delta_i = 4$. 
Let $\delta_{m_1}=\delta_{m_2}=\delta_{m_3}=1$ with 
$1 \leq m_1 < m_2 < m_3 \leq d$. 
Let $p_1=m_1-1$, $p_2=m_2-m_1-1$, 
$p_3=m_3-m_2-1$ and $p_4=d-m_3$. 
By $\delta_{i+1}=\delta_{d-i}$ for $0 \leq i \leq [(d-1)/2]$, 
one has $p_1=p_4$ and $p_2=p_3$. 
Moreover, by (\ref{Stanley}), one has $p_1 \geq p_2.$ 
Thus, 
\begin{eqnarray}\label{eqeq}
p_1 \geq p_2 \geq 0, \;\;\; 2p_1+2p_2=d-3. 
\end{eqnarray}
Our work is to construct an integral convex polytope $\Pc$ with dimension $d$
whose $\delta$-vector coincides with 
$$
\d(\Pc)=(1,\underbrace{0,\ldots,0}_{p_1},1,\underbrace{0,\ldots,0}_{p_2},
1,\underbrace{0,\ldots,0}_{p_2},1,\underbrace{0,\ldots,0}_{p_1})
$$
for an arbitrary integer $1 \leq m_1 < m_2 < m_3 \leq d$ 
satisfying the conditions (\ref{eqeq}). 
When $p_1=p_2=0$, it is easy to construct it 
by Examples \ref{example} (a). 
When $p_1=p_2>0$, if we set $d=4p_1+3$ and $m=4$, 
then the $\delta$-vector of the integral convex polytope 
whose vertices are of the form (\ref{vertex}) coincides with 
$\d(\Pc)=(1,\underbrace{0,\ldots,0}_{p_1},1,\underbrace{0,\ldots,0}_{p_1},
1,\underbrace{0,\ldots,0}_{p_1},1,\underbrace{0,\ldots,0}_{p_1})$
by virtue of Corollary \ref{corollary}. 

\begin{Lemma}
\label{first}
Let $d=4k+3$, $l \geq 1$ and $d'=d+2l$. 
There exists an integral simplex $\Pc \subset \RR^{d'}$
of dimension $d'$ whose $\delta$-vector coincides with
\begin{eqnarray}\label{d1}
(1,\underbrace{0,\ldots,0}_{k+l},1,\underbrace{0,\ldots,0}_{k},
1,\underbrace{0,\ldots,0}_{k},1,\underbrace{0,\ldots,0}_{k+l}) 
\in \ZZ^{d'+1}.
\end{eqnarray}
\end{Lemma}
\begin{proof}
When $l=1$, if we set $d=4k+5$ and $m=4$, 
then the $\delta$-vector of the integral convex polytope 
whose vertices are of the form (\ref{vertex}) 
coincides with (\ref{d1}). 
When $l \geq 2$, let $v_0',v_1',\ldots,v_{4k+2l+3}' \in \RR^{4k+2l+3}$ 
be the vertices as follows: 
\begin{eqnarray*}
v_i'=
\begin{cases}
(v_1,\underbrace{1,1,\ldots,1,1}_{2l-2}), \quad &i=1, \\
(v_i,\underbrace{0,1,\ldots,0,1}_{2l-2}), &i=2,3, \\
(v_i,\underbrace{0,0,\ldots,0,0}_{2l-2}), &i=4,\ldots,4k+5, \\
{\bf e}_i, &i=4k+6,\ldots,4k+2l+3, \\
(0,\ldots,0), &i=0, 
\end{cases}
\end{eqnarray*}
where $v_1,\ldots,v_{4k+5}$ are of the form (\ref{vertex}) 
with $d=4k+5$ and $m=4$. 
Then a simple computation enables us to show that 
\begin{eqnarray*}
\begin{vmatrix}
v_1' \\
\vdots \\
v_{4k+2l+3}'
\end{vmatrix}
=
\begin{vmatrix}
& &                & & &               & \\
& &\text{\Huge{A}} & & &\text{\huge{*}}& \\
& &                & & &               & \\
& &                & &1&               & \\
& &\text{\huge{0}} & & &\ddots         & \\
& &                & & &               &1
\end{vmatrix}
=4, 
\end{eqnarray*}
where $A$ is the determinant (\ref{matrix}) 
with $d=4k+5$ and $m=4$. 
One has 
\begin{eqnarray*}
\frac{3}{4}(v_0',1)+\frac{1}{4}\{(v_1',1)+\cdots+(v_{4k+5}',1)\}
+\frac{3}{4}\{(v_{4k+6}',1)+(v_{4k+8}',1)+\cdots+(v_{4k+2l+2}',1)\} \\
+\frac{1}{4}\{(v_{4k+7}',1)+(v_{4k+9}',1)+\cdots+(v_{4k+2l+3}',1)\} 
= (1,1,\ldots,1,k+l+1)
\end{eqnarray*}
and 
\begin{eqnarray*}
\frac{1}{2}\{(v_0',1)+(v_1',1)+\cdots+(v_{4k+2l+3}',1)\}
= (\underbrace{2,\ldots,2}_{4k+5},\underbrace{1,2,\ldots,1,2}_{2l-2},2k+l+2). 
\end{eqnarray*}
Hence $\delta_{k+l+1}=\delta_{2k+l+2}=\delta_{3k+l+3}=1$, 
as required. 
\end{proof}

Next, we consider the case of $\sum_{i=0}^{d} \delta_i = 5$. 
Let $\delta_{m_1}=\cdots=\delta_{m_4}=1$ with 
$1 \leq m_1 < \cdots < m_4 \leq d$. 
Let $p_1=m_1-1$, $p_2=m_2-m_1-1$, 
$p_3=m_3-m_2-1$,$p_4=m_4-m_3-1$ and $p_5=d-m_4$. 
By $\delta_{i+1}=\delta_{d-i}$ for $0 \leq i \leq [(d-1)/2]$, 
one has $p_1=p_5$ and $p_2=p_4$. 
Moreover, by (\ref{Stanley}), one has $p_1 \geq p_2,p_3.$ 
Thus, 
\begin{eqnarray}\label{eqeqeq}
p_1 \geq p_2,p_3 \geq 0, \;\;\; 2p_1+2p_2+p_3=d-4. 
\end{eqnarray}
Our work is to construct an integral convex polytope $\Pc$ with dimension $d$
whose $\delta$-vector coincides with 
$$
\d(\Pc)=(1,\underbrace{0,\ldots,0}_{p_1},1,\underbrace{0,\ldots,0}_{p_2},
1,\underbrace{0,\ldots,0}_{p_3},1,\underbrace{0,\ldots,0}_{p_2},1,\underbrace{0,\ldots,0}_{p_1})
$$
for an arbitrary integer $1 \leq m_1 < \cdots < m_4 \leq d$ 
satisfying the conditions (\ref{eqeqeq}). 
When $p_1=p_2=p_3=0$, it is easy to construct it 
by Examples \ref{example} (a). 
When $p_1=p_2=p_3>0$, if we set $d=5p_1+4$ and $m=5$, 
then the $\delta$-vector of the integral convex polytope 
whose vertices are (\ref{vertex}) coincides with 
$\d(\Pc)=(1,\underbrace{0,\ldots,0}_{p_1},1,\underbrace{0,\ldots,0}_{p_1},
1,\underbrace{0,\ldots,0}_{p_1},1,\underbrace{0,\ldots,0}_{p_1},1,\underbrace{0,\ldots,0}_{p_1})$
by virtue of Corollary \ref{corollary}. 

\begin{Lemma}
\label{second}
Let $d=5k+4$ and $l > 0$. \\
{\em (a)} Let $d'=d+2l$. 
There exists an integral simplex $\Pc \subset \RR^{d'}$
of dimension $d'$ whose $\delta$-vector coincides with
\begin{eqnarray}\label{d2}
(1,\underbrace{0,\ldots,0}_{k+l},1,\underbrace{0,\ldots,0}_{k},
1,\underbrace{0,\ldots,0}_{k},1,\underbrace{0,\ldots,0}_{k},
1,\underbrace{0,\ldots,0}_{k+l}) 
\in \ZZ^{d'+1}.
\end{eqnarray}
{\em (b)} Let $d'=d+3l$. 
There exists an integral simplex $\Pc \subset \RR^{d'}$
of dimension $d'$ whose $\d$-vector coincides with
\begin{eqnarray}\label{d3}
(1,\underbrace{0,\ldots,0}_{k+l},1,\underbrace{0,\ldots,0}_{k},
1,\underbrace{0,\ldots,0}_{k+l},1,\underbrace{0,\ldots,0}_{k},
1,\underbrace{0,\ldots,0}_{k+l}) 
\in \ZZ^{d'+1}.
\end{eqnarray}
{\em (c)} Let $d'=d+4l$. 
There exists an integral simplex $\Pc \subset \RR^{d'}$
of dimension $d'$ whose $\d$-vector coincides with
\begin{eqnarray}\label{d4}
(1,\underbrace{0,\ldots,0}_{k+l},1,\underbrace{0,\ldots,0}_{k+l},
1,\underbrace{0,\ldots,0}_{k},1,\underbrace{0,\ldots,0}_{k+l},
1,\underbrace{0,\ldots,0}_{k+l}) 
\in \ZZ^{d'+1}.
\end{eqnarray}
\end{Lemma}
\begin{proof}
A proof can be done as the similar way of Lemma \ref{first}. 

(a) When $l=1$, if we set $d=5k+6$ and $m=5$, 
then the $\delta$-vector of the integral convex polytope 
whose vertices are of the form (\ref{vertex}) 
coincides with (\ref{d2}). 
When $l \geq 2$, let $v_0',v_1',\ldots,v_{5k+2l+4}' \in \RR^{5k+2l+4}$ 
be the vertices as follows: 
\begin{eqnarray*}
v_i'=
\begin{cases}
(v_1,\underbrace{1,1,\ldots,1,1}_{2l-2}), \quad &i=1, \\
(v_i,\underbrace{0,1,\ldots,0,1}_{2l-2}), &i=2,3,4, \\
(v_i,\underbrace{0,0,\ldots,0,0}_{2l-2}), &i=5,\ldots,5k+6, \\
{\bf e}_i, &i=5k+7,\ldots,5k+2l+4, \\
(0,\ldots,0), &i=0, 
\end{cases}
\end{eqnarray*}
where $v_1,\ldots,v_{5k+6}$ are of the form (\ref{vertex}) 
with $d=5k+6$ and $m=5$. 
Then a simple computation enables us to show that 
$\begin{vmatrix}
v_1' \\
\vdots \\
v_{5k+2l+4}'
\end{vmatrix}
=5. $
One has 
\begin{eqnarray*}
\frac{4}{5}(v_0',1)+\frac{1}{5}\{(v_1',1)+\cdots+(v_{5k+6}',1)\}
+\frac{4}{5}\{(v_{5k+7}',1)+(v_{5k+9}',1)+\cdots+(v_{5k+2l+3}',1)\} \\
+\frac{1}{5}\{(v_{5k+8}',1)+(v_{5k+10}',1)+\cdots+(v_{5k+2l+4}',1)\} 
= (1,1,\ldots,1,k+l+1)
\end{eqnarray*}
and 
\begin{eqnarray*}
\frac{3}{5}(v_0',1)+\frac{2}{5}\{(v_1',1)+\cdots+(v_{5k+6}',1)\}
+\frac{3}{5}\{(v_{5k+7}',1)+(v_{5k+9}',1)+\cdots+(v_{5k+2l+3}',1)\} \\
+\frac{2}{5}\{(v_{5k+8}',1)+(v_{5k+10}',1)+\cdots+(v_{5k+2l+4}',1)\} 
= (\underbrace{2,\ldots,2}_{5k+6},\underbrace{1,2,\ldots,1,2}_{2l-2},2k+l+2). 
\end{eqnarray*}
Hence $\delta_{k+l+1}=\delta_{2k+l+2}=\delta_{3k+l+3}=\delta_{4k+l+4}=1$, 
as required. 

(b) When $l=1$, if we set $d=5k+7$ and $m=5$, 
then the $\delta$-vector of the integral convex polytope 
whose vertices are of the form (\ref{vertex}) 
coincides with (\ref{d3}). 
When $l \geq 2$, let $v_0',v_1',\ldots,v_{5k+3l+4}' \in \RR^{5k+3l+4}$ 
be the vertices as follows: 
\begin{eqnarray*}
v_i'=
\begin{cases}
(v_i,\underbrace{1,1,1,\ldots,1,1,1}_{3l-3}), \quad &i=1,2, \\
(v_i,\underbrace{0,1,1,\ldots,0,1,1}_{3l-3}), &i=3,4, \\
(v_i,\underbrace{0,0,0,\ldots,0,0,0}_{3l-3}), &i=5,\ldots,5k+7, \\
{\bf e}_i, &i=5k+8,\ldots,5k+3l+4, \\
(0,\ldots,0), &i=0, 
\end{cases}
\end{eqnarray*}
where $v_1,\ldots,v_{5k+7}$ are of the form (\ref{vertex}) 
with $d=5k+7$ and $m=5$. 
Then a simple computation enables us to show that 
$\begin{vmatrix}
v_1' \\
\vdots \\
v_{5k+3l+4}'
\end{vmatrix}
=5. $
One has 
\begin{eqnarray*}
\frac{3}{5}(v_0',1)+\frac{1}{5}\{(v_1',1)+\cdots+(v_{5k+7}',1)\}
+\frac{3}{5}\{(v_{5k+8}',1)+(v_{5k+11}',1)+\cdots+(v_{5k+3l+2}',1)\} \\
+\frac{1}{5}\{(v_{5k+9}',1)+(v_{5k+10}',1)+\cdots+(v_{5k+3l+3}',1)+(v_{5k+3l+4}',1)\} 
= (1,1,\ldots,1,k+l+1)
\end{eqnarray*}
and 
\begin{eqnarray*}
\frac{1}{5}(v_0',1)+\frac{2}{5}\{(v_1',1)+\cdots+(v_{5k+7}',1)\}
+\frac{1}{5}\{(v_{5k+8}',1)+(v_{5k+11}',1)+\cdots+(v_{5k+3l+2}',1)\} \\
+\frac{2}{5}\{(v_{5k+9}',1)+(v_{5k+10}',1)+\cdots+(v_{5k+2l+3}',1)+(v_{5k+2l+4}',1)\} \\
= (\underbrace{2,\ldots,2}_{5k+7},\underbrace{1,2,2,\ldots,1,2,2}_{3l-3},2k+l+2). 
\end{eqnarray*}
Hence $\delta_{k+l+1}=\delta_{2k+l+2}=\delta_{3k+2l+3}=\delta_{4k+2l+4}=1$, 
as required. 

(c) When $l=1$, if we set $d=5k+8$ and $m=5$, 
then the $\delta$-vector of the integral convex polytope 
whose vertices are of the form (\ref{vertex}) 
coincides with (\ref{d4}). 
When $l \geq 2$, let $v_0',v_1',\ldots,v_{5k+4l+4}' \in \RR^{5k+4l+4}$ 
be the vertices as follows: 
\begin{eqnarray*}
v_i'=
\begin{cases}
(v_i,\underbrace{1,1,1,1,\ldots,1,1,1,1}_{4l-4}), \quad &i=1,2,3, \\
(v_4,\underbrace{0,1,1,1,\ldots,0,1,1,1}_{4l-4}), &i=4, \\
(v_i,\underbrace{0,0,0,0,\ldots,0,0,0,0}_{4l-4}), &i=5,\ldots,5k+8, \\
{\bf e}_i, &i=5k+9,\ldots,5k+4l+4, \\
(0,\ldots,0), &i=0, 
\end{cases}
\end{eqnarray*}
where $v_1,\ldots,v_{5k+8}$ are of the form (\ref{vertex}) 
with $d=5k+8$ and $m=5$. 
Then a simple computation enables us to show that 
$\begin{vmatrix}
v_1' \\
\vdots \\
v_{5k+4l+4}'
\end{vmatrix}
=5. $
One has 
\begin{eqnarray*}
\frac{2}{5}(v_0',1)+\frac{1}{5}\{(v_1',1)+\cdots+(v_{5k+8}',1)\}
+\frac{2}{5}\{(v_{5k+9}',1)+(v_{5k+13}',1)+\cdots+(v_{5k+4l+1}',1)\} \\
+\frac{1}{5}\{(v_{5k+10}',1)+(v_{5k+11}',1)+(v_{5k+12}',1)+\cdots
+(v_{5k+4l+2}',1)+(v_{5k+4l+3}',1)+(v_{5k+4l+4}',1)\} \\
= (1,1,\ldots,1,k+l+1)
\end{eqnarray*}
and 
\begin{eqnarray*}
\frac{4}{5}(v_0',1)+\frac{2}{5}\{(v_1',1)+\cdots+(v_{5k+8}',1)\}
+\frac{4}{5}\{(v_{5k+9}',1)+(v_{5k+13}',1)+\cdots+(v_{5k+4l+1}',1)\} \\
+\frac{2}{5}\{(v_{5k+10}',1)+(v_{5k+11}',1)+(v_{5k+12}',1)+\cdots
+(v_{5k+4l+2}',1)+(v_{5k+4l+3}',1)+(v_{5k+4l+4}',1)\} \\
= (2,2,\ldots,2,2k+2l+2). 
\end{eqnarray*}
Hence $\delta_{k+l+1}=\delta_{2k+2l+2}=\delta_{3k+2l+3}=\delta_{4k+3l+4}=1$, 
as required. 
\end{proof}

{\bf Acknowledgements}

The author would like to thank Prof. T. Hibi 
for helping me in writing this paper.

\end{document}